\documentclass[a4paper,twoside,11pt]{article}
\usepackage{amssymb}
\usepackage{amsmath}
\usepackage{amsthm}
\usepackage{latexsym}
\usepackage{amsfonts}
\usepackage{graphicx}
\usepackage{graphics}

\theoremstyle{plain}
\newtheorem{thm}{Theorem}[section]   % Αρίθμηση συνεχόμενη (όχι κατά Θεώρημα, Λήμμα κ.λπ.)
\newtheorem{prop}[thm]{Proposition}
\newtheorem{lem}[thm]{Lemma}
\newtheorem{clm}[thm]{Claim}
\newtheorem{cor}[thm]{Corollary}
\newtheorem{Def}[thm]{Definition}%[section]

\theoremstyle{definition}
\newtheorem{rem}[thm]{Remark}

\newtheorem*{Proof}{Proof}

\newcommand{\el}{\ell}

\newcommand{\ra}{\;\rightarrow\;}

\newcommand{\de}{\delta }

\newcommand{\e}{\varepsilon }
\newcommand{\f}{\varphi}
\newcommand{\Fi}{\varPhi}

\newcommand{\vPsi}{\varPsi}

\newcommand{\la}{\lambda }

\newcommand{\si}{\sigma }

\newcommand{\ti}{\tau }

\newcommand{\oo}{\omega}
\newcommand{\C}{\mathbb{C}}
\newcommand{\E}{\mathbb{E}}

\newcommand{\R}{\mathbb{R}}

\newcommand{\N}{\mathbb{N}}

\newcommand{\bP}{\mathbb{P}}

\newcommand{\ssum}{\sum\limits}

\newcommand{\od}{\overline{d}}

\newcommand{\cm}{{\cal{M}}}
\newcommand{\cb}{{\cal{B}}}

\newcommand{\ld}{\ldots}

 \newcommand{\ud}{\underline{d}}

\begin{document}
\pagestyle{myheadings}
\markboth{Algebraic genericity of frequently universal harmonic functions on trees}{N. Biehler, V. Nestoridis, A. Stavrianidi}
\title{\bf Algebraic genericity of frequently universal harmonic functions on trees}
%
%\markright{Jordan domains with a rectifiable arc in their boundary}
%\titlerunning{Total unboundedness}
%
\author{N. Biehler, V. Nestoridis, A. Stavrianidi}
\date{}
\maketitle
\begin{abstract}
We show that the set of frequently universal harmonic functions on a tree $T$ contains a vector space except 0 which is dense in the space of harmonic functions on $T$ seen as subset of $\C^T$. In order to prove this we replace the complex plane $\C$ by any separable Fr\'{e}chet space $E$ and we repeat all the theory.
\end{abstract}
{\em AMS classification numbers}: 05C05, 31A20, 60J49, 30K99, 46M99 \smallskip
\section{Introduction}\label{sec1}
\noindent

For the theory of universality and the use of Baire's Category theorem we refer to \cite{9}, \cite{8}, \cite{13}, \cite{6}, \cite{12}, \cite{4}, \cite{7}, \cite{5}.

A tree $T$ is a connected, simply connected countable graph without non-trivial loops. With abuse of notation we shall also write $T$ for the set of vertices of the tree. We do not assume that $T$ is homogeneous: The number of edges joining at every vertex may be different, but shall always be finite. For $x,y\in T$ we write $x\sim y$ if $x$ and $y$ are neighbors, and also for any $x,y\in T$ there exists a unique natural number $n\in\N$ and a unique minimal finite sequence $\{z_0,z_1,\ld,z_n\}$ such that $z_0=x$, $z_n=y$ and $z_k\sim z_{k+1}$ for every $k<n$; this sequence is called the geodesic path from $x$ to $y$ and will be denoted by $[x,y]$.

The natural number $n$ as above is called the length of $[x,y]$ and will be denoted by $\el(x,y)$; $\el$ is a metric on $T$. By fixing a reference vertex $x_0\in T$ we obtain a partial ordering $\le$ as such: $x\le y$ if $x$ belongs to the geodesic path from $x_0$ to $y$. We denote by $y^-$ the unique neighbor of $y\neq x_0$ that belongs to $[x_0,y]$ that is $y^-\le y$. We call $y^-$ the parent of $y$. We shall denote with $T_n$ the set of all vertices at distance $n$ from $x_0$, and with $T_0=\{x_0\}$. For every $x\in T$ we also denote with $S(x)$ the set of all $y\in T$ such that $y\sim x$ and $y\neq y^-$; we will assume that for every vertex $x$, the set $S(x)$ contains at least two elements. For every $x\in T$ and every $y\in S(x)\cup\{x^-\}$ we associate a number $q(x,y)\ge0$ such that $q(x,x^-)=0$ and $q(x,y)>0$ for all $y\in S(x)$. These numbers satisfy:
\[
\sum_{y\in S(x)}q(x,y)=1.
\]
The set $T$ is infinite denumerable. The space $\C^T$ of all complex valued functions defined on $T$ is a complete metric space with respect to the metric:
\[
\rho(f,g)=\sum^\infty_{j=0}\frac{1}{2^j}\frac{|f(x_j)-g(x_j)|}{1+|f(x_j)-g(x_j)|}  \eqno{(\ast)}
\]
where $\{x_j:j\in\N\}$ is any enumeration of $T$. The induced topology is the topology of pointwise convergence.
\begin{Def}\label{def1.1} A function $f:T\ra\C$ is said to be $Q$-harmonic (or simply harmonic) if for all vertices $x\in T$ the following holds:
\[
f(x)=\sum_{y\sim x}q(x,y)f(y).
\]
\end{Def}

The set $H_Q(T)$ of all $Q$-harmonic functions on $T$ is a closed subset of $\C^T$ and thus, a complete metric space. (Therefore Baire's Category Theorem applies, \cite{9}, \cite{8}).

We define the boundary $\partial T$  of the tree $T$ to be the set of all the infinite geodesics starting at $x_0$; that is,  all $e=\{z_n:n\in\N\}$ such that $z_0=x_0$ and $z_{n+1}\in S(z_n)$ for every $n\in\N$. A vertex $x\in T$ defines a boundary sector $B_x$ of $\partial(T)$ which is the set of all $e\in\partial T$ such that $x\in e$. We assign a probability $p(B_x)$, where $x\in T_n$, as the product $\prod\limits^{n-1}_{j=0}q(y_j,y_{j+1})$, where $y_j\in T_j$, $y_n=x$, $y_0=x_0$ and of course $y_j\sim y_{j+1}$, for every $j$: $0\le j\le n-1$. The family $\{B_x\}_{x\in T_n}$ generates a $\si$-algebra $\cm_n$ on $\partial T$ and the above probabilities $p(B_x)$ extend to a probability measure $\bP_n$ on $\cm_n$. One can check that $\cm_n\subseteq\cm_{n+1}$ and $\bP_{n+1}|_{M_n}=\bP_n$ for every $n\in\N$. The Kolmogorov Consistency Theorem ([10]) implies that there exists a unique probability measure $\bP$ on the $\si$-algebra $\cm$ of subsets of $\partial T$ generated by the union $\bigcup\limits^\infty_{n=0}\cm_n$, such that $\bP|_{\cm_n}=\bP_n$ for every $n\in\N$.

Let $f\in H_Q(T)$ be a harmonic function. Then a sequence $\big\{\oo_n(f)\big\}^\infty_{n=0}$ of functions $\oo_n(f):\partial T\ra\C$ is defined by setting $\oo_n(f)(e)=f(y_n)$, where $y_n\in e\cap T_n$, for every $e\in\partial T$ and every $n\in\N$. It is obvious that the restriction of $\oo_n(f)$ on $\cb_x(x\in T_n)$ is constant and therefore, $\oo_n(f)$ is $\cm_n$-measurable and hence $\cm$-measurable.\\
Since $f(x)=\sum\limits_{y\sim x}q(x,y)f(y)$ it follows that the sequence $\big\{(\oo_n(f),\cm_n,\bP)\big\}^\infty_{n=0}$ defines a martingale and thus
\[
\E[\oo_m(f)\mid\cm_n]=\oo_n(f) \ \ \text{for every} \ \ m>n \ \ \text{\cite{11}}.
\]
\begin{Def}\label{def1.2}
Let $\ti\subseteq\N$ be an infinite set. A $Q$-harmonic function $f\in H_Q(T)$ belongs to the class $\widetilde{U}_\ti(T)$ if for every $\cm$-measurable function $\varPsi:\partial T\ra\C$ there exists a sequence $\big\{\la_n\big\}^\infty_{n=0}\subseteq\ti$ such that $\oo_{\la_n}(f)$ converges to $\varPsi$ in probability.
\end{Def}
\begin{rem}\label{rem1.3}
In Definition \ref{def1.2} it is equivalent to require that a subsequence $\oo_{\la_{k_n}}(f)$ of $\oo_{\la_n}(f)$ converge to $\varPsi$ $\bP$-almost surely. The space of $\cm$-measurable functions on $\partial T$ shall be denoted by $L^0(\partial T,\C)=L^0(\partial T)$ and considered with the topology induced by convergence in probability is a complete metric space with respect to the metric $d(\varPsi,\Fi)=\int\limits_{\partial T}\dfrac{|\varPsi-\Fi|}{1+|\varPsi-\Fi|}d\bP$. Definition \ref{def1.2} is analogue to the definition of universal trigonometric series in the sence of Menchoff (\cite{12}).
\end{rem}

In \cite{3} the following theorems are proven:
\begin{thm}\label{thm1.4}
The class $\widetilde{U}_\ti(T)$ is dense and $G_\de$ in $H_Q(T)$.
\end{thm}
\begin{thm}\label{thm1.5}
The class $\widetilde{U}_\ti(T)\cup\{0\}$ contains a dense linear subspace of $H_Q(T)$.
\end{thm}
\begin{Def}\label{def1.6}
For a subset $A$ of $\{0,1,2,\ld\}$ the lower density $\ud(A)$ of $A$ is defined as:
\[
\ud(A)=\liminf\frac{1}{n+1}\cdot card(\{m\in A:m\le n\}).
\]
\end{Def}

In \cite{3} the following definition has been given:
\begin{Def}\label{def1.7}
A harmonic function $f\in H_Q(T)$ is frequently universal, if for every non-void open subset $V$ of $L^0(\partial T)$, the set $\{n\in\{0,1,2,\ld\}:\oo_n(f)\in V\}$ has positive lower density. The set of frequently universal harmonic functions will be denoted as $U_{FM}(T)$.
\end{Def}

The following theorem is stated in \cite{3} without proof.
\begin{thm}\label{thm1.8}
The class $U_{FM}(T)$ is dense in $H_Q(T)$ and meager, that is, it is contained in a countable union of closed subsets with empty interior.
\end{thm}

The motivation of the present paper was to strengthen Theorem \ref{thm1.5} by proving that $U_{FM}\cup\{0\}$ contains a vector space which is dense in $H_Q(T)$. Initially we proved it by considering functions $f=(f_1,f_2,\ld)$, $f_k\in H_Q(T)$ harmonic on the tree $T$, $f:T\ra\C^\N$. The space of those $f$ is denoted by $(H_Q)^{N_0}$ or $H_Q(T,\C^{N_0})=(H_Q)^{N_0}$ endowed with the Cartesian topology. We use the notation: $\oo_n(f)=(\oo_n(f_1),\oo_n(f_2),\ld)$. Then $f$ is universal, if the sequence $\oo_n(f)$, $n=0,1,\ld$ is dense in the space $L^0(\partial T,\C^\N)$ of measurable functions $h:\partial T\ra\C^\N$, $h=(h_1,h_2,\ld)$,
\[
h_k\in L^0(\partial T,\C), \ \ L^0(\partial T,\C^\N)=[L^0(\partial T,\C)]^\N
\]
endowed with the Cartesian topology. Also $f$ is called frequently universal, if, for every non-empty open set $V\subset L^0(\partial T,\C^\N)$ the set:
$\{m\in\N:\oo_m(f)\in V\}$ has strictly positive lower density. If $f=(f_1,f_2,\ld)$ is frequently universal then it can be easily seen that the linear span of $f_1,f_2,\ld$ except zero, i.e. $span\{ f_1,f_2,\ld\}\backslash\{0\}$ is contained in the set of frequently universal functions $U_{FM}$. Furthermore, the sequence $f_1,f_2,\ld$ can be chosen to be dense in $H_Q(T)$ and simultaneously $f=(f_1,f_2,\ld)$ to be frequently universal. It follows that $span\{ f_1,f_2,\ld\}$ is a dense vector subspace of $H_Q(T)$ contained in $U_{FM}\cup\{0\}$. Thus, we have algebraic genericity for $U_{FM}$ which strengthens Theorem \ref{thm1.5}. It is also possible to repeat all the theory for the new harmonic universal functions $f:T\ra\C^\N$ $f=(f_1,f_2,\ld)\in(H_Q)^{N_0}=H_Q(T,\C^{N_0})$.

In order to prove algebraic genericity of the frequently universal functions $f=(f_1,f_2,\ld)\in(H_Q)^{N_0}=H_Q(T,\C^{N_0})$ we are led to consider the space $\big((H_Q)^{N_0}\big)^{N_0}=H_Q(T,(\C^{N_0})^{N_0})=(f_{n,k})_{n,k=1,2,\ld}$ and to repeat the theory for this new space and so on.

The spaces $\big(((H_Q^{N_0}))^{N_0}\ld\big)^{N_0}$ are separable Fr\'{e}chet spaces and everything will be simplified if we fix $E$ a separable Fr\'{e}chet space (on the field of real numbers $\R$ or complex numbers $\C$), if we consider the space $H_Q(T,E)=\{f:T\ra E$ harmonic$\}\subset E^T$ endowed with the cartesian topology and if we define what it means that $f\in H_Q(T,E)$ is universal or frequently universal and repeat all the theory. This is done in the following sections and as a corollary we obtain the algebraic genericity of $U_{FM}$.
\section{Universal harmonic functions}\label{sec2}
\noindent

Let $T$ be a tree as in \S 1 and let $E$ be a separable Fr\'{e}chet space. We denote by $d$ the metric in $E$.
\begin{Def}\label{def2.1}
A function $f:T\ra E$ is called harmonic if
\[
f(x)=\sum_{y\in S(x)}q(x,y)f(y) \ \ \text{for each} \ \ x\in T.
\]
The set of such functions is denoted by $H_Q(T,E)$. It is a subset of $E^T$ endowed with the cartesian topology.
\end{Def}

Since $E$ is complete and $T$ denumerable, those spaces are complete.

Since $E$ is a separable metric space, the same is true for $E^T$ and $H_Q(T,E)$.

If $f\in H_Q(T,E)$, the function $\oo_n(f):\partial T\ra E$ is defined in a similar way as in \S 1 and the sequence $\oo_n(f)$, $n=0,1,2,\ld$ is a martingale with values in $E$.

A function $h:\partial T\ra E$ is $\cm$-measurable ($\cm_n$-measurable) if, for every open set $W\subset E$ we have $h^{-1}(W)\in \cm$ $(h^{-1}(W)\in\cm_n$ respectively).

Identifying two $M$-measurable functions if they are $P$-almost everywhere equal, we obtain the space $L^0(\partial T,E)$.

We endow $L^0(\partial T,E)$ with the topology of convergence in probability, that is, for any $h_1,h_2\in L^0(\partial T,E)$ their distance is
\[
\widetilde{\rho}(h_1,h_2)=\int_{\partial T}\frac{d(h_1(x),h_2(x))}{1+d(h_1(x),h_2(x))}dP(x).
\]

Since $E$ is separable, for every $n$ the space of $\cm_n$-measurable functions defined on $\partial T$ with values in $E$ is separable. It follows that $L^0(\partial T,E)$ is also separable and it has a dense sequence $h_n\in L^0(\partial T,E)$, $n=1,2,\ld$ where each $h_n$ is $\cm_{\ti(n)}$-measurable for some $\ti(n)\in\{0,1,2,\ld\}$.
\begin{Def}\label{def2.2}
A harmonic function $f\in H_Q(T,E)$ is universal if the sequence $\oo_n(f)$, $n=0,1,2,\ld$ is dense in $L^0(\partial T,E)$. The set of universal harmonic functions $f\in H_Q(T,E)$ is denoted by $U(T,E)\subset H_Q(T,E)$.
\end{Def}
\begin{thm}\label{thm2.3}
Under the above assumptions and notation the set $U(T,E)$ is a dense and $G_\de$ subset of $H_Q(T,E)$; hence $U(T,E)$ is non-void.
\end{thm}
\begin{Proof}
The proof is similar to the proof of Theorem 2.2 in \cite{3} and is omitted. It suffices to replace $\C$ by $E$.
\end{Proof}
\begin{rem}\label{rem2.4}
Let $f\in U(T,E)$; then, the sequence $h_n=\oo_n(f)$, $n=0,1,2,\ld$ is dense in $L^0(\partial T,E)$ and $h_n$ is $\cm_n$-measurable.
\end{rem}
\begin{rem}\label{rem2.5}
Let $E,F$ be two isomorphic separable Fr\'{e}chet spaces and $\f:E\ra F$ an isomorphism. Let $f:T\ra E$ be a function. Then, $f\in H_Q(T,E)$ iff $\f\circ f\in H_Q(T,F)$ and $f\in U(T,E)$ iff $\f\circ f\in U(T,F)$.
\end{rem}
\section{Frequently universal harmonic functions}\label{sec3}
\noindent

In analogy with Definition \ref{def1.7} we give the following definition.
\begin{Def}\label{def3.1}
A function $f\in H_Q(T,E)$ is frequently universal, if, for every non-empty open set $V\subset L^0(\partial T,E)$ the set $\{n\in\N:\oo_n(f)\in V\}$ has strictly positive lower density. We denote the set of such functions by $U_{FM}(T,E)$. Obviously $U_{FM}(T,E)\subseteq U(T,E)$.
\end{Def}
\begin{rem}\label{rem3.2}
If $E$ and $F$ are two isomorphic separable Fr\'{e}chet spaces and $\f:E\ra F$ is an isomorphism, then $f\in U_{FM}(T,E)$ iff $\f\circ f\in U_{FM}(T,F)$.
\end{rem}
\begin{thm}\label{thm3.3}
Under the above assumptions and notation $U_{FM}(T,E)$ is dense and meager in $H_Q(T,E)$. In particular $U_{FM}(T,E)\neq\emptyset$.
\end{thm}

The proof of Theorem \ref{thm3.3} is a modification of the proof in \cite{1}, \cite{2}. Since these papers are not yet published we give all the details of the proof.\vspace*{0.2cm}\\
\noindent
{\bf Proof of Theorem \ref{thm3.3}}. We consider a dense sequence $\{h_k,\;k\in\N\}$ in $H_Q(T,E)$ where each $h_k$ is $\cm_k$-measurable. It can be easily seen that such a sequence exists according to Remark \ref{rem2.4}.

To construct a frequently universal function $f$ it suffices that, for every $k\in\N$ the set: $\{n\in\N:\oo_n(f)\in B\Big(h_k,\dfrac{1}{2^k}\Big)\Big\}$ has positive lower density. For this purpose we will use the following sequence of integers.
\begin{lem}\label{lem3.4}
Let $\el(k)=s+1$, where $k=2^s\cdot d$ and $d$ is an odd number. We also set $r_1=0+\el(1)$ and for $k\ge2$, $r_k=r_{k-1}+\el(k)$.

This sequence has the following properties:
\begin{enumerate}
\item[(i)] $r_n=\ssum^n_{k=1}\el(k)$
\item[(ii)] $r_{2^N}=2^{N+1}-1$ for $N\in\N$
\item[(iii)] For $m\le N$ it holds that $card(\{k:1\le k\le2^N$ such that $\el(k)=m\})=2^{N-m}$
\item[(iv)] For $m\ge1$ the set $\{r_n:n\in\N,\; \el(n)=m\}$ has strictly positive lower density.
\end{enumerate}
\end{lem}
\noindent
{\bf Proof of Lemma \ref{lem3.4}}. (\cite1),\cite{2}).

i) This is obvious.

ii) To see this we proceed by induction on $\N$: For $N=0,1$ it is obviously true. Suppose that $r_{2^N}=2^{N+1}-1$, we will prove that $r_{2^{N+1}}=2^{N+2}-1$.

We have that
\[
r_{2^{N+1}}=\sum^{2^{N+1}}_{k=1}\el(k)=\sum^{2^N-1}_{k=1}\el(k)+\el(2^N)+
\sum^{2^{N+1}-1}_{k=2^N+1}\el(k)+\el(2^{N+1}).
\]
But $\el(2^N)=N+1$, $\el(2^{N+1})=N+2$. We set
\[
S_1=\sum^{2^N-1}_{k=1}\el(k) \ \ \text{and} \ \ S_2=\sum^{2^{N+1}-1}_{k=2^N+1}\el(k).
\]
We have that $S_1=S_2$ because for $k=1,2,\ld,2^N-1$ we have $\el(k)=\el(k+2^N)$. Indeed let $k=2^s\cdot d$ with $s<N$ and $d$ odd. Then $\el(k)=s+1$ and $k+2^N=2^s\cdot d+2^N=2^s(d+2^{N-s})$ with $d+2^{N-s}$ an odd number since $N>s$. Thus $\el(k+2^N)=s+1=\el(k)$. It follows that:
\[
r{_{2^{N+1}}}=2\bigg(\sum^{2^N-1}_{k=1}\el(k)\bigg)+2(N+1)+1=2\bigg(\sum^{2^N}_{k=1}
\el(k)\bigg)+1
\]
because $N+1=\el(2^N)$.

We also have $\ssum^{2^N}_{k=1}\el(k)=r_{2^N}$ and therefore: $r_{2^{N+1}}=2(r_{2^N})+1$ and using the induction hypothesis we find $r_{2^{N+1}}=2(2^{N+1}-1)+1=2^{N+2}-1$.

This proves property (ii)

iii) For $m=1$ half of the $k\in[1,2^N]$, $k\in\N$ satisfy $\el(k)=1$ (the odd $k$), so
\[
card(\{k:1\le k\le 2^N:\el(k)=1\})=\frac{2^N}{2}=2^{N-1}.
\]
For $m=2$ we have half of the half (that is $\dfrac{1}{4})$ which give $\dfrac{2^N}{4}=2^{N-2}$.

Continuing this way we see that:
\[
card(\{k:1\le k\le 2^N:\el(k)=m\})=2^{N-m} \ \ \text{where} \ \ m\le N.
\]

iv) We want to show that the set:
\[
\{r_n:n\in\N;\el(n)=m\} \ \ \text{has strictly positive lower density in $\N$}.
\]
The sequence $r_n$ is increasing so we have to show that
\[
\liminf_{M\in\N}\frac{card\{n\in\N:r_n\le M,\el(n)=m\}}{M+1}>0
\]
Because $r_{2^N}=2^{N+1}-1$ it suffices to show the following:
\[
\inf_{N:N\ge m}\frac{card\{n:r_n\le r_{2^N},\el(n)=m\}}{r^{2^N}+1}>0.
\]
Since the sequence $r_n$ is strictly increasing this reduces to the following
\[
\inf_{N:N\ge m}\frac{card\{n:1\le n\le 2^N:\el(n)=m\}}{r^{2^N}+1}>0.
\]
Indeed, we have that the numerator is equal to $2^{N-m}$ as we have already seen. The denominator equals $2^{N+1}$. Thus, we find $\dfrac{2^{N-m}}{2^{N+1}}=2^{-m-1}>0$.

Thus, the set $\{r_n:n\in\N$ with $\el(n)=m\}$ has strictly positive lower density.

The proof of Lemma \ref{lem3.4} is completed.\medskip

To construct a function $f\in U_{FM}(T,E)$ it suffices to construct $f\in H_Q(T,E)$ such that:
\[
\oo_{r_k}(f)\in B\bigg(h_{\el(k)},\frac{1}{2^{\el(k)}}\bigg) \ \ \text{for all} \ \ k\ge1,
\]
where the sequence $h_n$, $n=0,1,\ld$ is dense in $H_Q(T,E)$ and each $h_n$ is $M_n$-measurable, as previously mentioned. For this purpose, we will need the following lemma.
\begin{lem}\label{lem3.5}
Let $n>0$, $s>0$ be natural numbers and let $\f$ be harmonic on $T_0\cup T_1\cup\cdots\cup T_s$. Then it is possible to extend $\f$ to a harmonic function $\vPsi$ on $T_0\cup T_1\cup\cdots\cup T_{s+n-1}$ defined on $T_0\cup T_1\cup\cdots\cup T_{s+n}$ (i.e. $\vPsi|_{T_0\cup T_1\cup\cdots\cup T_s}\equiv \f)$ in such a way that:
\[
\oo_{s+n}(\vPsi)\in B\bigg(h_n,\frac{1}{2^n}\bigg).
\]
\end{lem}
\noindent
{\bf Proof of Lemma \ref{lem3.5}}. (\cite{1}, \cite{2}).
Consider a vertex $v\in T_s$. Let $V\subseteq T_{s+n}$ be the set of all vertices of $T_{s+n}$ which are at distance $n$ from $v$. Clearly $V$ contains a vertex $w$ such that the probability to go from $v$ to $w$ is less than or equal to $\dfrac{1}{2^n}$.

We set
\[
\vPsi(u)=h_n(e) \ \ \text{where} \ \ e\cap T_{s+n}=\{u\}
\]
for every $u\in V\backslash\{w\}$. $\vPsi$ is well-defined this way because $h_n$ is $M_{n+s}$-measurable and hence the value of $h_n(e)$, where $e\cap T_{s+n}=\{u\}$ does not depend on the choice of $e$.

Clearly we can extend $\vPsi$ to all vertices $x\in T_{s+1}\cup\cdots\cup T_{s+n}$ such that $v\in[0,x]$ in such a way that $\vPsi$ is harmonic.

We do the same procedure with all vertices of $T_s$ and we define $\vPsi$ on $T_0\cup\cdots\cup T_{s+n}$ which is harmonic on $T_0\cup\cdots\cup T_{s+n-1}$ and extends $\f$. We also have: $\oo_{s+n}(\vPsi)\in B\Big(h_n,\dfrac{1}{2^n}\Big)$ because the set $\{w=w(v)\in T_{n+s}:v\in T_s\}$ has probability at most $\dfrac{1}{2^n}$ and outside this set we have $\vPsi\equiv h_n$ on $T_{n+s}$. This proves Lemma \ref{lem3.5}.\medskip

By Lemma \ref{lem3.5} it is possible to construct $f\in H_Q(T,E)$ by induction such that
\[
\oo_{r_k}(f)\in B\bigg(h_{\el(k)},\frac{1}{2^{\el(k)}}\bigg) \ \ \text{for every} \ \ k\in\N.
\]
Indeed, it suffices to apply Lemma \ref{lem3.5} for $s=r_{k-1}$, $n=\el(k)$, $s+n=r_k$ because $r_k=r_{k-1}+\el(k)$.

We have that $\oo_{r_k}(f)\in B\Big(h_{\el(k)},\dfrac{1}{2^{\el(k)}}\Big)$ for each $k\in\N$. Since, according to Lemma \ref{lem3.4}, the set $\{r_k,k\in\N:\el(k)=n\}$ has positive lower density, $f$ is frequently universal.

Hence $U_{FM}(T,E)\neq\emptyset$. We will now prove that $U_{FM}(T,E)$ is dense.

We consider $f_0\in U_{FM}(T,E)$. We also consider $\f\in H_Q(T,E)$ such that $\oo_n(\f)=\oo_{n_0}(\f)$ for every $n\ge n_0$ for some $n_0\in\N$. We can easily see that $f_0+\f\in U_{FM}(T,E)$ by considering the open set $V-\oo_{n_0}(\f)$ and using the fact that the set $\{n\in\N:\oo_n(f_0)\in V-\oo_{n_0}(\f)\}$ has positive lower density. Since the set $L=\{\f\in H_Q(T,E):\oo_n(\f)=\oo_{n_0}(\f)$ for all $n\ge n_0$ for some $n_0(\f)\in\N\}$ is dense in $H_Q(T,E)$ the same holds for the set of translations $f_0+\f$, $\f\in L$. This set is dense and contained in $U_{FM}(T,E)$. This proves the density of $U_{FM}(T,E)$.

We will now prove that $U_{FM}(T,E)$ is meager in $H_Q(T,E)$. We need the following preliminaries:
\begin{Def}\label{def3.6}
For a subset $A$ of the natural numbers $\N$ the upper density $\od(A)$ of $A$ is defined as:
\[
\od(A)=\limsup_{n\in\N}\frac{1}{n+1}\cdot card(\{m\in A:m\le n\}).
\]
\end{Def}
\begin{Def}\label{def3.7}
A harmonic function $f\in H_Q(T,E)$ belongs to the class $X(T,E)$ if the following holds:\\
For every non-void open subset $V$ of $L^0(\partial T,E)$ the set $\{n\in\N:\oo_n(f)\in V\}$ has upper density in $\N$ exactly equal to 1.
\end{Def}
\begin{prop}\label{prop3.8}
The class $X(T,E)$ is $G\de$ and dense in $H_Q(T,E)$.
\end{prop}
\noindent
{\bf Proof of Proposition \ref{prop3.8}}. (\cite{1}, \cite{2}).

Let $V_j$, $j\in\N$ be a denumarable base of open sets in $L^0(\partial T,E)$ and $n$, $m\in\N$. We consider the set:
$E(j,m,n)=\{f\in H_Q(T,E)$: There exists $q\in\N:n>q>\Big(1-\dfrac{1}{m}\Big)\cdot n$ and natural numbers $n_1,n_2,\ld,n_q$ such that $1\le n_1<n_2<\cdots<n_q<n$: $\oo_{n_\el}(f)\in V_j$ for all $\el=1,2,\ld,q\}$.\\
Then it can be seen that:
\[
X(T,E)=\bigcap_{j,m,N}\bigcup_{n\ge N}E(j,m,n).
\]
We can easily see that $E(j,m,n)$ is open for every $j,m,n\in\N$.

Therefore $X(T,E)$ is $G\de$ in $H_Q(T,E)$. In order to use Baire's Category Theorem in the complete metric space $H_Q(T,E)$ it suffices to prove that $\bigcup\limits_{n\ge N}E(j,m,n)$ is dense in $H_Q(T,E)$. We imitate the proof of the density in Theorem \ref{thm2.3}.

We construct a harmonic function $f$ up to a certain level $n_1$ such that $\oo_{n_1}(f)\in V_j$. We extend harmonically $f$ so that $\oo_k(f)=\oo_{n_1}(f)$ for all $k=n_1,n_1+1,\ld,n_1+\si$ where $\si>\Big(1-\dfrac{1}{m}\Big)(n_1+\si)$. To do this extension we set $f(x)=f(x^-)$ for all $x\in\bigcup\limits_{\el\ge n_1+1}T_\el$. This extension is harmonic because $q(x,x^-)=0$.

This proves Proposition \ref{prop3.8}.
\begin{prop}\label{prop3.9}
\[
X(T,E)\cap U_{FM}(T,E)=\emptyset.
\]
\end{prop}
\noindent
{\bf Proof of Proposition \ref{prop3.9}}. (\cite{1}, \cite{2}).
We suppose that $X(T,E)\cap U_{FM}(T,E)\neq\emptyset$ to arrive at a contradiction.

Let $f$ be an element of $X(T,E)\cap U_{FM}(T,E)$. We consider $V_1,V_2$ two disjoint open balls in $L^0(\partial T,E)$.\\
Then $\{n\in\N:\oo_n(f)\in V_1\}\cap\{n\in\N: \oo_n(f)\in V_2\}=\emptyset$.

Since $f\in U_{FM}(T,E)$, there exists $\de>0$ and $m\in\N$, such that $card(\{n\in\N:n\le N:\oo_n(f)\in V_1\})>\dfrac{\de}{2}(N+1)$ for all $N\ge m$.

Since $f\in X(T,E)$ there exists $N_0\ge m$ such that
\[
card(\{n\le N_0:\oo_n(f)\in V_2\})>\bigg(1-\frac{\de}{2}\bigg)(N_0+1).
\]
Combining those two we find:
\begin{align*}
1+N_0=&card(\{n\le N_0\})\ge card(\{n\le N_0:\oo_n(f)\in V_1\cup V_2\}) \\
=&card(\{n\le N_0:\oo_n(f)\in V_1\}) \\
&+card(\{n\le N_0:\oo_n(f)\in V_2\}) \\
>&\frac{\de}{2}(N_0+1)+\bigg(1-\frac{\de}{2}\bigg)(N_0+1)=N_0+1.
\end{align*}
Thus: $1+N_0>1+N_0$ which is absurd.

This proves Proposition \ref{prop3.9}.\medskip

Since $U_{FM}(T,E)$ is disjoined from $X(T,E)$, which is $G\de$ and dense, we conclude that $U_{FM}(T,E)$ is meager in $H_Q(T,E)$. This completes the proof of Theorem \ref{thm3.3}.
\section{Algebraic genericity}\label{sec5}
\noindent

Let $E$ be a separable Fr\'{e}chet space, then the infinite denumerable product $E^\N$ is also a separable Fr\'{e}chet space. According to Theorem \ref{thm3.3} there exists $f=(f_1,f_2,\ld)\in U_{FM}(T,E^\N)$ where $f_k:T\ra E$ are harmonic functions.
\begin{clm}\label{clm4.1}
The linear space $\langle f_1,f_2,\ld\rangle$ of $f_1,f_2,\ld$ where $f=(f_1,f_2,\ld)\in U_{FM}(T,E^\N)$ is contained in $U_{FM}(T,E)\cup\{0\}$.
\end{clm}
\begin{Proof}
We consider a linear combination of those $f_n$, $n\in\N$, $a_1f_1+\cdots+a_sf_s$ where $a_s\neq0$.

Let $V$ be a non-void, open subset of $L^0(\partial T,E)$. There exists some function $h$ such that $B\Big(h,\dfrac{1}{2^M}\Big)\subseteq V$ for some $M\in\N$.

Let $\e=\dfrac{1}{s\cdot 2^M}$ and we consider the set:
\begin{align*}
\widetilde{V}=&B\bigg(0,\e\cdot\min\bigg\{\frac{1}{|a_1|},1\bigg\}\bigg)\times
B\bigg(0,\e\cdot\min\bigg\{\frac{1}{|a_2|},1\bigg\}\bigg) \\
&\times\cdot\times B\bigg(\frac{h}{|a_s|},\e\cdot\min\bigg\{\frac{1}{|a_s|},1\bigg\}\bigg) \\
&\times L^0(\partial T,E)\times \ld\;.
\end{align*}
Since the metric $(\widetilde{\rho})$ in $L^0(\partial T,E^\N)$ is equivalent to the cartesian metric considering that $L^0(\partial T,E^\N)$ is the infinite denumerable product of spaces $L^0(\partial T,E)$, we have that $\widetilde{V}$ is open in $L^0(\partial T,E^\N)$.

(We assume $\dfrac{1}{0}=\infty$ and hence $\min\{\infty,1\}=1$).

Since $f\in U_{FM}(T,E^\N)$ we have that the set $\{n\in\N:\oo_n(f)\in\widetilde{V}\}$ has positive lower density. Therefore, for all $n$ in this set we have
\[
\oo_n(f)\in\widetilde{V}\Leftrightarrow\oo_n(f_i)\in B\bigg(0,\e\cdot\min\bigg\{\frac{1}{|a_i|},1\bigg\}\bigg)
\]
for $1\le i\le s-1$ and
\[
\oo_n(f_s)\in B\bigg(\frac{h}{|a_s|},\e\cdot\min\bigg\{\frac{1}{|a_s|},1\bigg\}\bigg).
\]
It follows that: $\oo_n(a_if_i)\in B(0,\e)$ for $1\le i\le s-1$ and $\oo_n(a_sf_s)\in B(h,\e)$.

As a result
\[
\widetilde{\rho}(\oo_n(a_if_i),0)<\e \ \ \text{for} \ \ 1\le i\le s-1 \ \ \text{and} \ \ \widetilde{\rho}(\oo_n(a_sf_s),h)<\e.
\]
We now have:
\begin{align*}
&\widetilde{\rho}(\oo_n(a_1f_1+\cdots+a_sf_s),h)\le\widetilde{\rho}(\oo_n(a_1f_1),0)+\cdots+\widetilde{\rho}
(\oo_n(a_{s-1}f_{s-1}),0) \\
&+\widetilde{\rho}(\oo_n(a_sf_s),h)\le s\cdot\e=\frac{1}{2^M}
\end{align*}
where we used the fact that the metric $d$ in the Fr\'{e}chet space $E$ is invariant under translations.

This implies $\oo_n(a_1f_1+\cdots+a_sf_s)\in B\Big(h,\dfrac{1}{2^M}\Big)\subseteq V$ and this is true for every $n\in\N$ such that $\oo_n(f)\in\widetilde{V}$. The set of such $n\in\N$ has positive lower density. As a result the set $\{n\in\N:\oo_n(a_1f_1+\cdots+a_sf_s)\in V\}$ has positive lower density. This implies that $a_1f_1+\cdots+a_sf_s\in U_{FM}(T,E)$ and thus
\[
span\{f_n:,n\in\N\}\subseteq U_{FM}(T,E)\cup\{0\}.
\]
The proof of the claim is completed.
\end{Proof}

Next we will show the following.
\begin{clm}\label{clm4.2} There exists a sequence $f_1,f_2,\ld$ dense in $H_Q(T,E)$ such that $f=(f_1,f_2,\ld)\in U_{FM}(T,E^\N)$.
\end{clm}
\noindent
{\bf Proof of the claim}. We consider a sequence $\{\f_n:n\in\N\}$ which is dense in $H_Q(T,E)$ and a function $f\in U_{FM}(T,E^\N)$, $f=(f_1,f_2,\ld)$.

For every $n\in\N$ there exists some $j_0(n)\in\N$ such that:
\[
\sum^\infty_{j=j_0(n)}\frac{1}{2^j}<\frac{1}{n}.
\]
For each $n\in\N$ we set $g_n(x_j)=\f_n(x_j)-f_n(x_j)$ for all $j=0,1,\ld,j_0(n)-1$.

There also exists some $N(n)$ such that $\{x_0,x_1,\ld,x_{j_0(n)}\}\subseteq\bigcup\limits^{N(n)}_{k=0}T_k$.

Then we also set $g_n(x)=\f_n(x)-f_n(x)$ for every $x\in\bigcup\limits^{N(n)}_{k=0}T_k$ and all $n\in\N$.

Finally, we set $g_n(x)=g_n(x^-)$ for every $x\in T_k$, $k>N(n)$. We notice that $\oo_k(g_n)=\oo_{N(n)}(g_n)$ for every $k>N(n)$.

As a result: $g_n+f_n\in B\Big(\f_n,\dfrac{1}{n}\Big)$ for each $n\in\N$.

Since the set $\{\f_n:n\in\N\}$ is dense in $H_Q(T,E)$ we deduce that the set of functions: $\{g_n+f_n:n\in\N\}$ is also dense in $H_Q(T,E)$.

We will now show that the function $f+g=(f_1+g_1,f_2+g_2,\ld)\in H_Q(T,E^\N)$ is frequently universal.

The metric $\widetilde{\rho}$ in the space $L^0(\partial T, E^\N)$ is equivalent to the cartesian metric when considering $L^0(\partial T, E^\N)$ as the infinite denumerable product space of $L^0(\partial T,E)$. We have that a basic open set $V$ in the space $(L^0(\partial T, E^\N),\widetilde{\rho})$ can be written as
\[
V=V_1\times V_2\times\cdots\times V_s\times L^0(\partial T,E)\times \ld
\]
where $V_i$ $i=1,\ld,s$ are non-void open sets in $L^0(\partial T,E)$.

We have that $\oo_n(g_1)=\oo_{N_1}(g_1)$ for $n\ge N_1,\ld,\oo_n(g_s)=\oo_{N(s)}(g_s)$ for $n\ge N_s$.

We set $L=\max\{N(j):1\le j\le s\}$. Then for all $n\ge L$ and $j=1,2,\ld s$ we have that $\oo_n(g_j)=\oo_L(g_j)$.

Since $f\in U_{FM}(T,E^\N)$ and $V-\oo_L(g)$ is open we have:
\[
\underline{d}(\{n\in\N:\oo_n(f)\in V-\oo_L(g)\})>0.
\]
Equivalently:
\[
\begin{array}{l}
        \underline{d}(\{n\in\N:\oo_n(f)+\oo_L(g)\in V\})>0 \ \ \text{or} \\       [1ex]
        \underline{d}(\{n\ge L:\oo_n(f)+\oo_L(g)\in V\})>0 \ \ \text{or}   \\ [1ex]
  \underline{d}(\{n\ge L:\oo_n(f)+\oo_n(g)\in V\})>0 \ \ \text{or}       \\ [1ex]
\underline{d}(\{n\ge L:\oo_n(f+g)\in V\})>0 \ \ \text{or}         \\[1ex]
        \underline{d}(\{n\in\N:\oo_n(f+g)\in V\})>0.
      \end{array}
\]
Since $V$ was an arbitrary basic open set in $L^0(\partial T,E^\N)$ it follows that $f+g\in U_{FM}(T,E^\N)$.

The proof of the claim is complete. \medskip

Combining Claims \ref{clm4.1} and \ref{clm4.2} we obtain
\begin{thm}\label{thm4.3}
Under the above assumptions and notation $U_{FM}(T,E)\cup\{0\}$ contains a vector space which is dense in $H_Q(T,E)$.
\end{thm}
\begin{Proof}
Let $f=(f_1,f_2,\ld)\in U_{FM}(T,E^\N)$ as in Claim \ref{clm4.2}. Since the sequence $f_1,f_2,\ld$ is dense in $H_Q(T,E)$ it follows that their linear span $\langle f_1,f_2,\ld\rangle$ is also dense in $H_Q(T,E)$. Since $f=(f_1,f_2,\ld)\in U_{FM}(T,E^\N)$, it follows that
\[
\langle f_1,f_2,\ld\rangle\subset U_{FM}(T,E)\cup\{0\}
\]
according to Claim \ref{clm4.1}. This completes the proof of Theorem \ref{thm4.3}.
\end{Proof}
\begin{cor}\label{cor4.4} Under the above assumptions and notation $U(T,E)\cup\{0\}$ contains a vector space which is dense in $H_Q(T,E)$.
\end{cor}
\begin{Proof}
This follows immediately from Theorem \ref{thm4.3} because $U_{FM}(T,E)\subset U(T,E)$ Corollary \ref{cor4.4} for $E=\C$ coincides with Theorem \ref{thm1.5}. Theorem \ref{thm4.3} for $E=\C$ gives a strengthened version of Theorem \ref{thm1.5} which is the algebraic genericity of $U_{FM}$ and was the motivation of this paper.

\end{Proof}
\bigskip
\newpage
\noindent
National and Kapodistrian University of Athens\\
Department of Mathematics\\
Panepistemiopolis, 157 84\\
Athens,
Greece \bigskip\\
e-mail addresses: \\
nikiforosbiehler@outlook.com \\
vnestor@math.uoa.gr \\
aleksandrastavrianidi@yahoo.gr

\end{document}